\documentclass[11pt]{amsart}
\usepackage[top=1.4in, bottom=1.4in, left=1in, right=1in]{geometry}
\usepackage{graphicx}
\usepackage{amsmath,latexsym}
\usepackage{amsfonts,amssymb}
\usepackage{amscd, amsthm}
\usepackage{stmaryrd}
\usepackage{fancyhdr}
\usepackage{commath,enumitem,chngcntr}
\usepackage[titletoc,toc,title]{appendix}
\usepackage{bbm}
\usepackage{pst-node}
\usepackage{tikz-cd} 
\usepackage{pgf,tikz}
\usepackage{mathrsfs}
\usetikzlibrary{arrows,calc,positioning}
\usepackage{float}

\makeatletter
\def\BState{\State\hskip-\ALG@thistlm}
\makeatother

\usepackage{pst-node}

\usepackage{tikz-cd} 
\usepackage{pgf,tikz}
\usepackage{mathrsfs}
\usetikzlibrary{arrows,calc,positioning}
\usepackage{float}
\usepackage{overpic}

 \newtheorem{theorem}{Theorem}[section]
\newtheorem{lemma}[theorem]{Lemma}
\newtheorem{proposition}[theorem]{Proposition}

\usepackage{lineno}

\title{A Note on Geometric Triangulations of the Hyperbolic Plane}


\begin{document}

\author{Yanwen Luo}
\address{Department of Mathematics, Oklahoma State University, Stillwater, OK, 74074}
\email{yanwen.luo@okstate.edu}
\author{Tianqi Wu}

\address{Clark University, Worcester, MA, 01610
}

\email{tianwu@clarku.edu}

\author{Xiaoping Zhu}

\address{Department of Mathematics, Rutgers University, Piscataway, NJ, 08854}

\email{xz349@scarletmail.rutgers.edu}



\begin{abstract}
It is well-known that the Euclidean plane has a standard 6-regular triangulation by equilateral triangles and the hyperbolic plane has a $k$-regular geometric triangulation for every $k>6$. 
In this note, we construct 6-regular geometric triangulations of the hyperbolic plane with a uniform bound on edge lengths, resolving a conjecture by Feng Luo \cite{Luo}.  We also show that there is no uniform bound on stretch factors of Delaunay triangulations in the hyperbolic plane, in contrast to the existence of such a  bound in the Euclidean plane.
\end{abstract}
\maketitle

\section{Introduction}
A triangulation $T$ on a surface $S$ is \textit{$k$-regular} for some integer $k>0$ if the degrees of all its vertices are $k$, and it is \textit{geometric} if all the edges are embedded as geodesic arcs with respect to a Riemannian metric on $S$. In this note, we consider the existence of $k$-regular geometric triangulations on the unit $2$-sphere $\mathbb{S}^2$, the Euclidean plane $\mathbb{E}^2$, and the hyperbolic plane $\mathbb{H}^2$.

The most familiar degree-regular geometric triangulation is the standard hexagonal triangulation in the Euclidean plane $\mathbb{E}^2$ tiled by equilateral triangles. It is straightforward to construct $k$-regular triangulations of $\mathbb{S}^2$ with $k = 3, 4, 5$ by projecting the regular tetrahedron, octahedron, and icosahedron to $\mathbb{S}^2$. In the hyperbolic plane $\mathbb{H}^2$, we can use triangle groups 
to produce $k$-regular geometric triangulations for any $k>6$.  
\begin{figure}[!ht]
\label{degreeregular}
\includegraphics[width=0.6\linewidth]{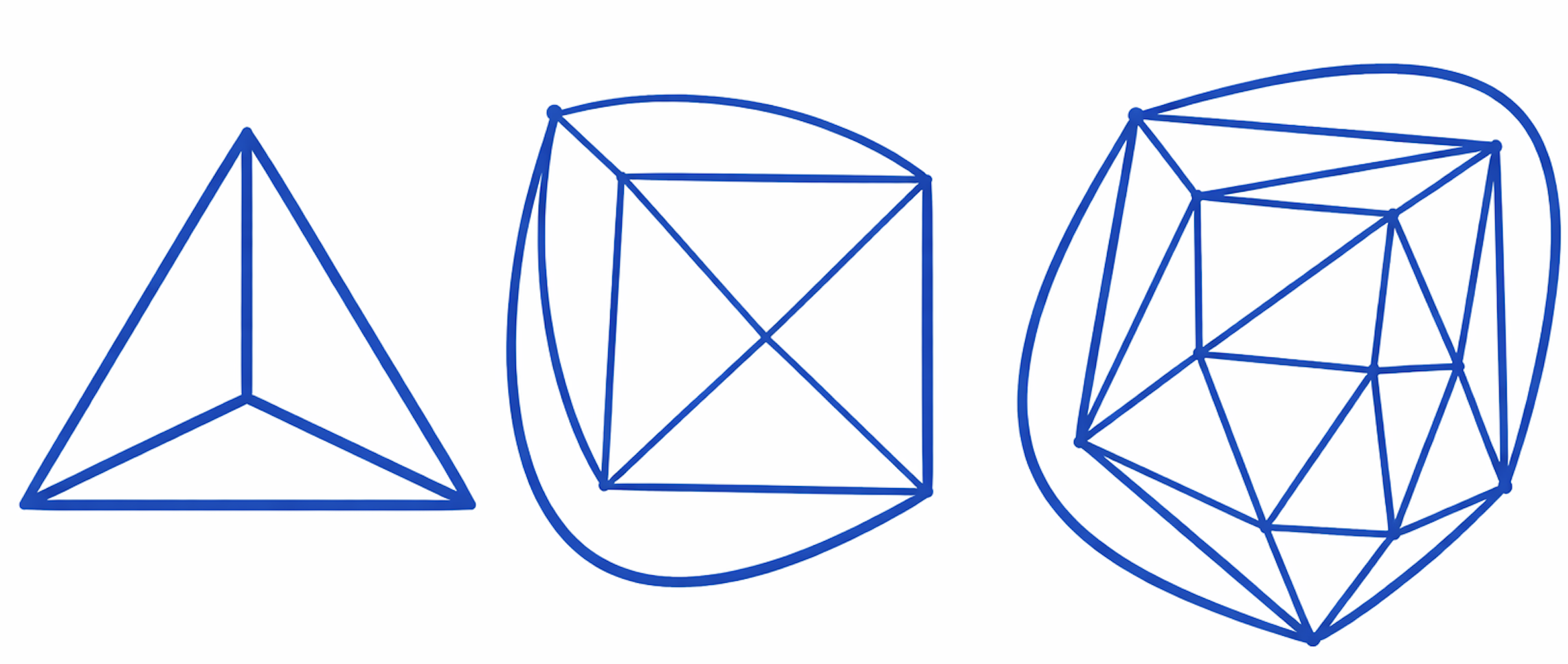}
  \includegraphics[width=0.25\linewidth]{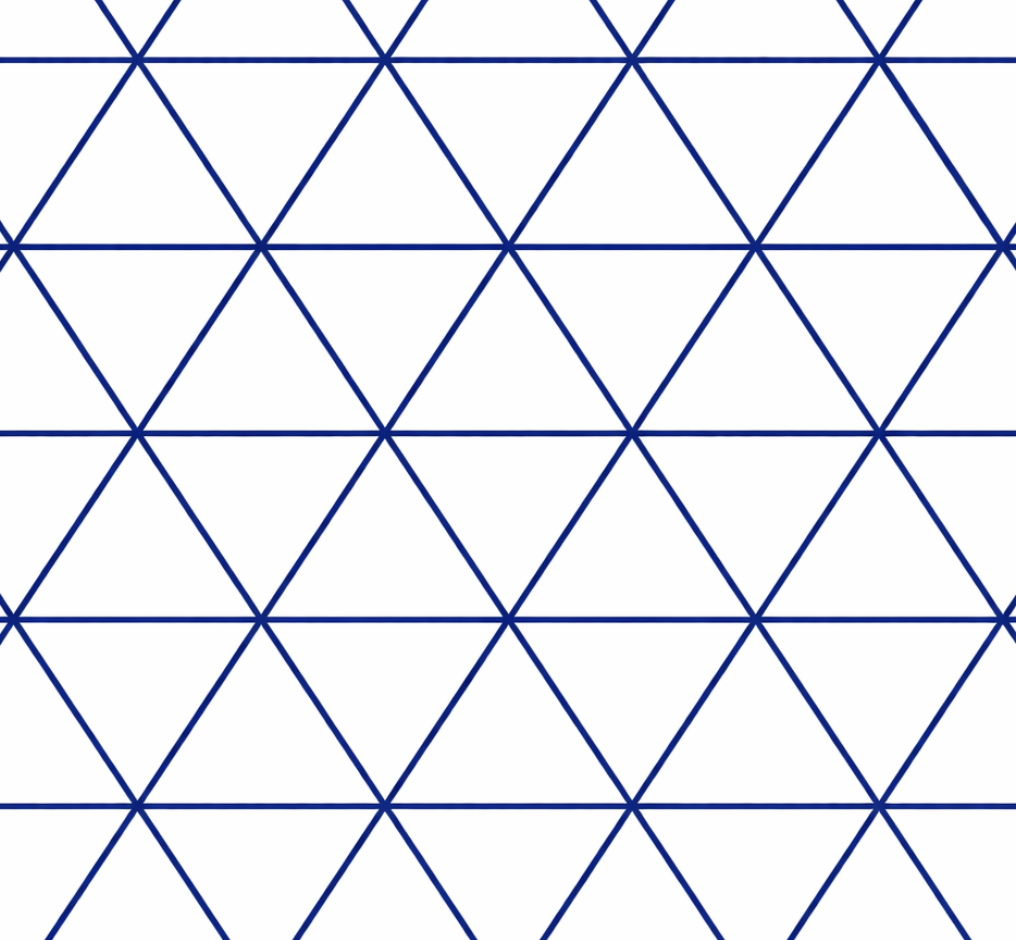}
  \caption{$k$-regular geometric triangulations with $k = 3, 4, 5, 6$.}
\end{figure}

For any $k\geq 6$, the Euler characteristic immediately implies that there is no topological $k$-regular triangulation on $\mathbb{S}^2$. It is a natural question whether the similar phenomenon holds for $\mathbb{E}^2$ and $\mathbb{H}^2$. In particular, it is asked by Luo \cite{Luo}  whether there is a $6$-regular geometric triangulation of $\mathbb{H}^2$ with a uniform bound on diameters, or equivalently edge lengths of triangles. 

Intuitively, such triangulations may not exist. If $k = 6$, the number of vertices grows linearly with respect to the combinatorial distance from one vertex, while it grows exponentially if $k>6$, reflecting different growth rates of the Euclidean geometry and hyperbolic geometry. However, 
\begin{theorem}
\label{main}
    There is a $6$-regular geometric triangulation in $\mathbb{H}^2$ such that diameters of its triangles are uniformly bounded. There is a $k$-regular geometric triangulation in $\mathbb{E}^2$ for $k\geq 7$ such that diameters of its triangles are uniformly bounded. 
\end{theorem}

The problem above arises from the dichotomy between the geometry of the Euclidean plane and the hyperbolic plane in various aspects. The celebrated \textit{uniformization theorem} of simply connected Riemann surfaces states that the Euclidean plane and the hyperbolic plane are not conformally equivalent, which gives rise to the concept of \textit{type} of Riemann surfaces.


Circle packings proposed by Thurston \cite{Thurston} as discrete conformal maps lead to a similar concept of the \textit{CP-type} of infinite triangulations of the plane: each triangulation of the topological plane corresponds to an infinite circle packing of either the Euclidean plane $\mathbb{E}^2$ or the hyperbolic plane $\mathbb{H}^2$, but not both \cite{HS1995}. In particular, there is no $6$-regular circle packing in $\mathbb{H}^2$. It is also related to the concepts of recurrence and transience in network theory. See an excellent survey \cite{Bowers}.

Moreover, hyperbolic spaces serve as fundamental models to classify metric spaces up to quasi-isometry, leading to the concept of hyperbolic groups in geometric group theory. Clearly, the Euclidean plane and the hyperbolic plane are not quasi-isometric. The problem above asks a combinatorial version of this fact. 
 Theorem \ref{main} shows that the condition on the $6$-regular triangulation is not sufficiently rigid. Hence, we may consider other more restrictive classes of triangulations. 

\textit{Delaunay triangulations} are promising candidates as fundamental objects in computational geometry \cite{De}. It has been shown in \cite{CRS} that every Delaunay triangulation of the Euclidean plane is not a hyperbolic metric space. 
Furthermore, the stretch factor of Delaunay triangulation in the Euclidean plane is uniformly bounded \cite{Xia}. These evidences suggest that Delaunay triangulations, as metric graphs, can be regarded as discrete approximations to the Euclidean plane up to quasi-isometry.

However, we construct explicit examples in the hyperbolic plane showing that Delaunay triangulations in the hyperbolic plane do not have such a property. 
\begin{theorem}
\label{main2}
    There is no uniform bound on  stretch factors of Delaunay triangulations in  $\mathbb{H}^2$. 
\end{theorem}
We will prove some basic facts in Section 2 and construct triangulations in Theorem \ref{main} in Section 3 and 4. The definition of stretch factors and the proof of Theorem \ref{main2} will be given in Section 5.


\section{Preliminary facts}
\subsection{Regular triangulations of small degrees.}
It is a simple observation that there is no $k$-regular triangulation on $\mathbb{E}^2$ and $\mathbb{H}^2$ if $k = 3, 4, 5$.

\begin{lemma}
There is no $k$-regular triangulation on $\mathbb{E}^2$ and $\mathbb{H}^2$ if $k = 3, 4, 5$.
\end{lemma}
\begin{proof}
     Let $v_0$ be a vertex in $T$ and $v_1,\cdots, v_k$ be its neighbors in its star, denoted by $star(v_0)$. Since $T$ is a triangulation, $v_1,\cdots, v_k$ with $v_1 = v_{k+1}$ are distinct vertices. If $k=3$, then we can not add more edges or vertices in $T$, since all the vertices have degree $3$ in $star(v_0)$. It is the graph of the tetrahedron. 
    
    If $k=4$, there is one edge in $T$ not in $star(v_0)$ from each of $v_1, v_2, v_3, v_4$ respectively, and there is a triangle in $T$ sharing the edge $v_iv_{i+1}$ for each $i = 1, 2, 3, 4$ containing this extra edge. Notice that this edge can not connect two vertices among these four. Indeed, say $v_1v_3$ is an edge in $T$, then one of $v_2$ and $v_4$ has degree $3$. There is a vertex $v_5$ in $T$ not in $star(v_0)$, and it is immediate that $v_1, \cdots, v_4$ must be connected to $v_5$, otherwise there is a vertex of degree larger than $4$. This leads to the graph of the octahedron.

    If $k=5$, by the same argument as that of $k=4$, there are five distinct triangles sharing the edge $v_iv_{i+1}$ with five triangles $v_0v_iv_{i+1}$, denoted by $v_iv_{i+1}u_i$. Notice that $u_1, u_2, \cdots, u_5$ are distinct, otherwise one of $v_1, \cdots, v_5$ has degree $4$. Notice that $u_1, u_2, \cdots, u_5$ have degree at least $4$ in the induced graph $G$ by $v_0, v_1, \cdots, v_5, u_1, \cdots u_5$. Then there is an extra edge from each of $u_1, \cdots, u_5$ not in $G$. The same above shows that these edges can not connect any pair of $u_1, \cdots, u_5$. Hence, there is one extra vertex $u_0$ and all $u_1, \cdots, u_5$ are connected to $u_0$, which produces the graph of the icosahedron. 
\end{proof}

\subsection{Regular geometric triangulations}
To construct $k$-regular geometric triangulations in $\mathbb{E}^2$ with $k>6$, and $6$-regular geometric triangulations in $\mathbb{H}^2$, we first point out that such triangulations contain slim triangles.

\begin{lemma}
Let $T$ be a $k$-regular geometric triangulation on $\mathbb{E}^2$ with $k>6$, or a $6$-regular geometric triangulation on $\mathbb{H}^2$. Then the infimum of angles in $T$ is zero. 
\end{lemma} 

\begin{proof}
	If not, in the first case for $\mathbb{E}^2$ we can construct a quasiconformal homeomorphism from $\mathbb{E}^2$ to $\mathbb{H}^2$ by sending a triangle in $T$ to the corresponding triangle in the standard geometric triangulation of $\mathbb{H}^2$ generated by reflecting equilateral triangles with inner angle $2\pi/k$. Similarly, in the second case we can construct a quasiconformal homeomorphism from $\mathbb{H}^2$ to $\mathbb{E}^2$ sending $T$ to the standard equilateral hexagonal triangulation. This contradicts the fact that there is no quasiconformal homeomorphism from $\mathbb{E}^2$ to $\mathbb{H}^2$.
\end{proof}


\section{$k$-regular geometric triangulations on the Euclidean plane}
It is known that one can construct $k$-regular geometric triangulations on $\mathbb{E}^2$ for $k>6$.
The idea is to construct chains of circles with radius $n$ and distribute the points on the circle evenly. The similar computation is given in \cite{DG2021} . We  show that these triangulations have the desired uniform bounds. 
\begin{lemma}
For any $k\geq 7$, there exist $k$-regular geometric triangulations on $\mathbb{E}^2$ with uniform upper bounds on the lengths of edges. 
\end{lemma}
\begin{proof}
	The construction of a $k$-regular geometric triangulation is given as follows.  Let $C_n$ be circles with radius $r_n = n$ centered at the origin of $\mathbb{E}^2$ and $T_k$ be the topological $k$-regular triangulation of $\mathbb{E}^2$. Let $a_n$ with $n\geq 1$ be the number of vertices whose combinatorial distance to a vertex $v_0$ in $T_k$ is exactly $n$. Notice that these vertices form a cycle in $T_k$. We call them \textit{vertices in the $n$th layer}. This cycle consists of $a_n$ edges, so we can regard $a_n$ as the number of edges in the $n$th layer of $T_k$. We set $a_0 = 0$ since there is no edge connecting $v_0$ with itself. We call edges connecting vertices in the same layer \textit{angular} and edges connecting vertices in different layers \textit{radial}.

We have the following recurrence relation
$$a_{n+1} = (k-4)a_n - a_{n-1}$$
with $a_0 = 0$ and $a_1 = k$. Indeed, consider the star of the $n$th layer vertices with $n\geq 2$. It is an annulus consisting of vertices from the $(n-1)$th, $n$th, and $(n+1)$th layers. Each triangle in this annulus contains exactly one angular edge and two radial edges. Each angular edge in the $(n-1)$th layer and the $(n+1)$th layer corresponds to exactly one triangle, and each angular edge in the $n$th layer corresponds to exactly two triangles. Hence, we have 
$$F = a_{n-1} + a_{n+1} + 2a_n,$$
where $F$ is the number of triangles in this annulus. Meanwhile, since every vertex in the $n$th layer has degree $k$, there are exactly $(k-2)a_n$ radial edges in this annulus. Therefore, 
$$2F = 2(k-2)a_n,$$
and we obtain the recurrence relation above. 

Solve this relation with the initial conditions, and we have the formula for $a_n$ with $n\geq 1$
$$a_n = \frac{k}{\sinh \lambda}\sinh(n\lambda), \quad \cosh \lambda = \frac{k-4}{2}.$$

    We put $v_0$ at the origin and the vertices in the $1$st layer to be a regular $k$-gon on $C_1$. Then we construct the $k$-regular geometric triangulation inductively. Assume that $v_n^i$ for $i = 1, \cdots, a_n$ with $v_n^1 = v_n^{a_n+1}$ are vertices in the $n$th layer placed in $C_{n}$. Then we determine the locations of $v_{n+1}^j$ for $j = 1, \cdots, a_{n+1}$, and $v_{n+1}^{a_{n+1}+1} = v_{n+1}^1$. 

    For each edge $v_n^iv_n^{i+1}$, there is a ray $\gamma$ from the origin passing through the midpoint of $v_n^iv_n^{i+1}$ intersecting with $C_{n+1}$ at a unique point $p_i$. This point will be a vertex in the $(n+1)$th layer. Then we can determine $a_n$ vertices among $a_{n+1}$ vertices of the $(n+1)$th layer, each of which determines a triangle $v_n^iv_n^{i+1}p_i$.

\begin{figure}[ht!]
\centering
\begin{overpic}[width=0.6\linewidth]{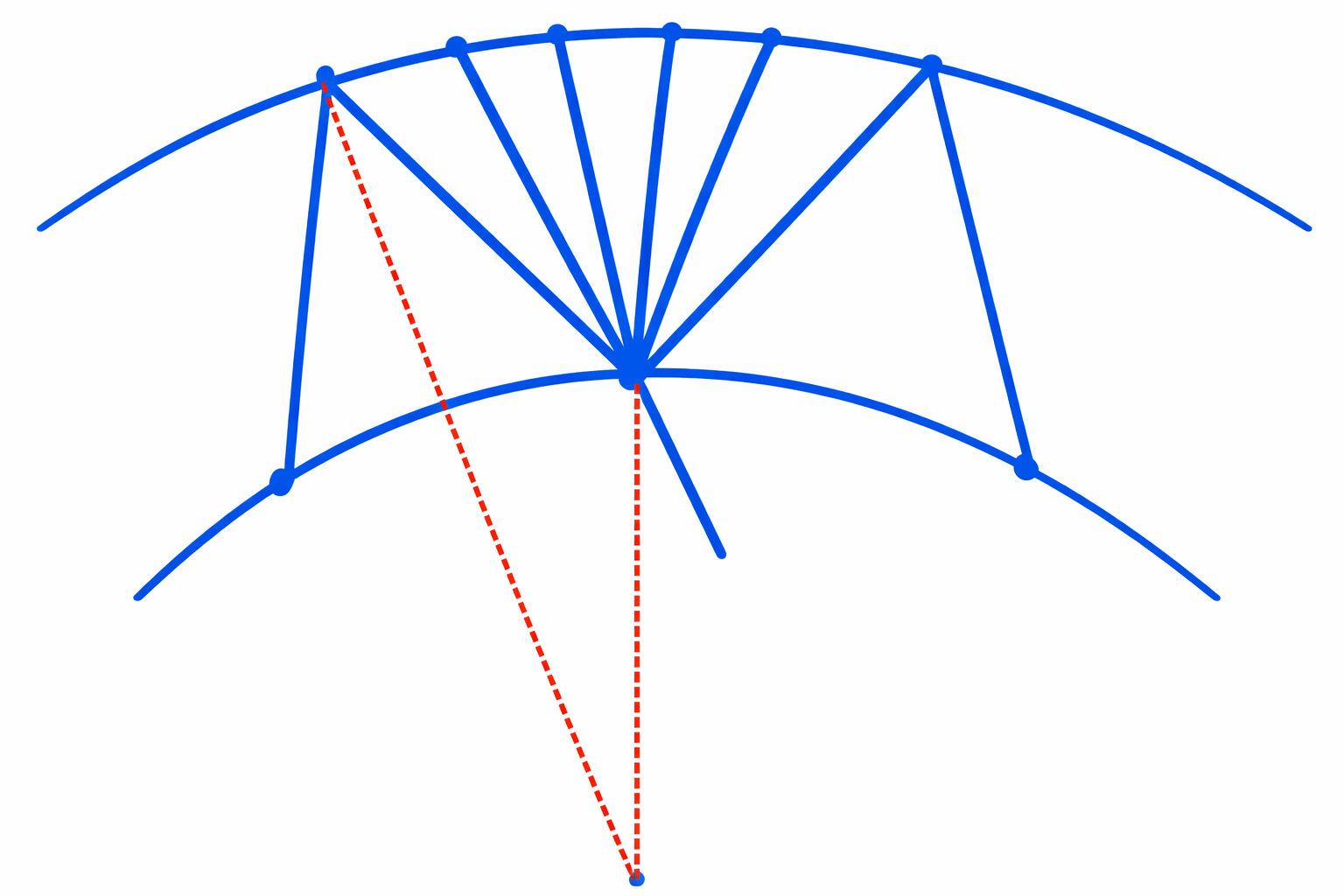}
    \put(50,0){$v_0$}
    \put(23,64){$p_{i-1}$}
    \put(20,27){$v^{i-1}_n$}
        \put(70,27){$v^{i+1}_n$}
        \put(42,34){$v^{i}_n$}
    \put(68,64){$p_{i}$}

\end{overpic}
 \caption{Radial edges between vertices in the $n$th and $(n+1)$th layer.}
  \label{einduction}
\end{figure}

    We  distribute the neighbors of $v_{n}^i$ in the $(n+1)$th layer in the arc $p_{i-1}p_i$ in $C_{n+1}$ evenly and connect the consecutive vertices in $C_{n+1}$. Then we can connect $v_n^i$ with the vertices in the arc $p_{i-1}p_i$ to determine all radial edges. Notice that this will not introduce intersections of edges, since the polygon determined by $v_0$, $p_i, p_{i-1}$ and vertices in the arc $p_{i-1}p_i$ is a convex polygon, and $v_n^i$ is in this convex polygon by the construction. 
    
    Also, notice that there are $(k-4)$ or $(k-3)$ vertices in $C_{n+1}$ connected with each vertex in $C_n$ by this construction. Hence, there is always some vertex between $p_i$ and $p_{i+1}$ in $C_{n+1}$, and 
    $$\theta_{n+1} = \max_{i=1, \cdots, a_{n+1}} \theta(v_{n+1}^iv_{n+1}^{i+1}) \leq \frac{1}{2}\max_{i=1, \cdots, a_{n}} \theta(v_{n}^iv_{n}^{i+1}) = \frac{\theta_n}{2}, $$
    where $\theta(v_{n}^iv_{n}^{i+1})$ is the  length of the arc between $v_{n}^i$ and $v_{n}^{i+1}$ in $C_n$. 
    Then $\theta_n = c2^{-n}$ for some constant $c>0$. The length of the circular arc between $v_{n}^i$ and $v_{n}^{i+1}$ is bounded by $r_n\theta_n\leq cn2^{-n}$. Hence, it tends to zero as $n\to \infty$, and the lengths of all angular edges tend to zero. 

    Meanwhile, the lengths of radial edges between $n$th and $(n+1)$th layer are bounded by lengths of isosceles triangles such as $p_iv_n^iv_n^{i+1}$. It can be estimated directly from the cosine law  
    $$l(p_iv_{n}^i)^2 \leq n^2 + (n+1)^2 - 2n(n+1)\cos \theta_n \leq n^2 + (n+1)^2 - 2n(n+1)\cos\frac{c}{2^n} \leq 1 + \frac{c^2n(n+1)}{2^n},$$
    which is clearly bounded by $2$ for sufficiently large $n$. We also see that the lengths of radial edges tend to $1$.

\end{proof}

\section{6-regular geometric triangulations on the hyperbolic plane}

We now give a construction of the desired $6$-regular geometric triangulations with uniform bounds in the hyperbolic plane using the Klein model of $\mathbb{H}^2$.
\begin{figure}[h!]
  \includegraphics[width=0.7\linewidth]{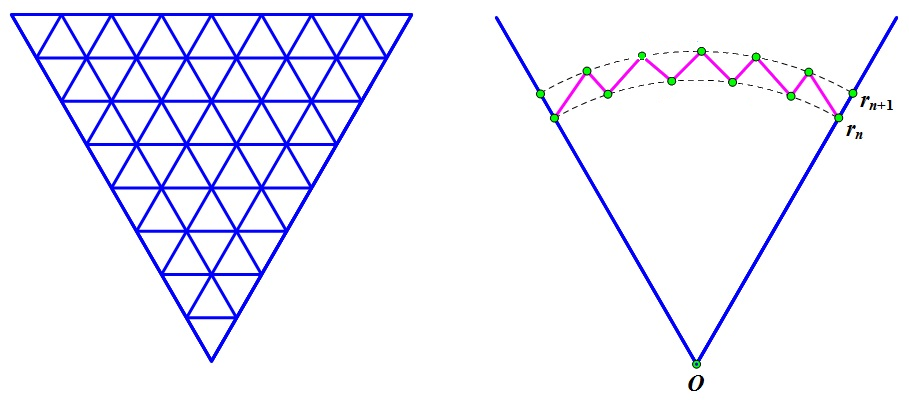}
  \caption{Edges of type 2 between vertices in the $n$th layer and the $(n+1)$th layer.}
\end{figure}

\begin{enumerate} 
	\item Put one vertex $v_0$ in the origin, and shoot out two rays $\gamma_1$ and $\gamma_2$ along the positive $x$-axis and along the direction $e^{i\pi/3}$.
	\item Let $r_n = \alpha \log (n+1)$ for $n\geq 1$ with $0<\alpha<1/2$. Construct a sequence of points $a_n$ and $b_n$ on $\gamma_1$ and $\gamma_2$ with distance $r_n$ from $v_0$ in the hyperbolic metric. Connect $a_n$ with $a_{n+1}$ and $b_n$ with $b_{n+1}$ using geodesics. Call these edges of \textit{type 0}. 
	\item  Connect $a_n$ with $b_n$ by circle arcs $C_n$ of radius $r_n$, called the $n$-th layer in this construction. Evenly distribute $n-1$ points on $C_n$. Then each $C_n$ is divided into $n$ arcs. Connect the consecutive vertices on $C_n$ to generate $n$ edges. Call these edges of \textit{type 1}. By symmetry, these edges have the same length.
	\item Connect points on $C_n$ and $C_{n+1}$ based on the combinatorics of the hexagonal triangulation to generate $2n$ edges. Call these edges of \textit{type 2}.
	\item Rotate the configuration by $k\pi/3$ for $k = 1, \cdots, 5$ and identify the edges of type 0 to generate the full hexagon geometric triangulation. 
\end{enumerate}
We can give explicit coordinates of the vertices on $C_n$ in this triangulation in the plane for $n\geq 1$
$$v_n^k = \tanh r_ne^{i\frac{k\pi}{3n}}, \quad k = 0, 1, \cdots, n.$$
And the vertex $v_n^k$ is connected with vertices $v_{n+1}^k$ and $v_{n+1}^{k+1}$ in the next layer. 
\begin{proposition}
\ The construction above generates a 6-regular geometric triangulation of $\mathbb{H}^2$ with uniform bound on the lengths of edges. 
\end{proposition}
The proof of this proposition consists of two lemmas. 

\begin{lemma}
The construction above generates a valid geometric triangulation with no intersection of edges and degenerate triangles. 
\end{lemma}
\begin{proof}
	The proof is based on induction on $n$. The base case $n = 1$ is trivial. Assume it is true for the triangulation generated up to the $n$-th layer. We will show that the triangles between $C_n$ and $C_{n+1}$ lie in the annulus  $A_n$ bounded  $C_n$ and $C_{n+1}$. There are two types of triangles:
	\begin{enumerate}
	
		\item Type 1: two vertices  on $C_n$ and one vertex on $C_{n+1}$. The three vertices are 
		$$v_n^k =\tanh r_ne^{i\frac{k\pi}{3n}}, v_n^{k+1} = \tanh r_ne^{i\frac{(k+1)\pi}{3n}}, \text{and }v_{n+1}^{k+1} =\tanh r_{n+1}e^{i\frac{(k+1)\pi}{3(n+1)}}.$$
		Notice that 
		$$\frac{k}{n} < \frac{k+1}{n+1} < \frac{k+1}{n}, \text{and } r_{n}\leq r_{n+1},$$
        then other than the edges of type 0, all the other edges are not colinear with the origin. And two edges of type 2 in the triangle will not intersect at the interior, since their interiors are contained in disjoint open sectors from the origin. 
		\item Type 2: one vertex on $C_n$ and two vertices on $C_{n+1}$. In this case, the three vertices are 
		$$v_n^k =\tanh r_ne^{i\frac{k\pi}{3n}}, v_{n+1}^{k+1} = \tanh r_{n+1}e^{i\frac{(k+1)\pi}{3(n+1)}}, \text{ and }v_{n+1}^{k} = \tanh r_{n+1}e^{i\frac{k\pi}{3(n+1)}}.$$
		Notice that 
		$$\frac{k}{n+1} < \frac{k}{n} < \frac{k+1}{n+1}, \text{ and } r_{n}\leq r_{n+1}.$$
        Similarly, these edges of type 2 will not intersect. Hence, any pair of edges of type 2 in $A_n$ are disjoint in their interior. Notice that they are also disjoint in the interiors from the edges of type 1 in the $(n+1)$th layer along $C_{n+1}$ by the convexity of the polygon determined by the origin and $v_{n+1}^k$ for $k = 0, 1, \cdots, n+1$.
	\end{enumerate}
	
Hence, we only need to show that the edges of type 2 lie in the ring $A_n$, so it does not intersect with previous layers. 
	
%
%
The critical figures for triangles of type 1 and type 2 are given below, where the angles at $v_n^k$ or $v_{n}^{k+1}$ are Euclidean right angles in the triangle $v_n^kv_{n+1}^kv_{n+1}^{k+1}$ (middle) and $v_n^{k+1}v_{0}v_{n+1}^{k+1}$ (left) respectively in Figure \ref{caseanalysis}. 

	\begin{figure}[H]
    \includegraphics[width=0.32\linewidth]{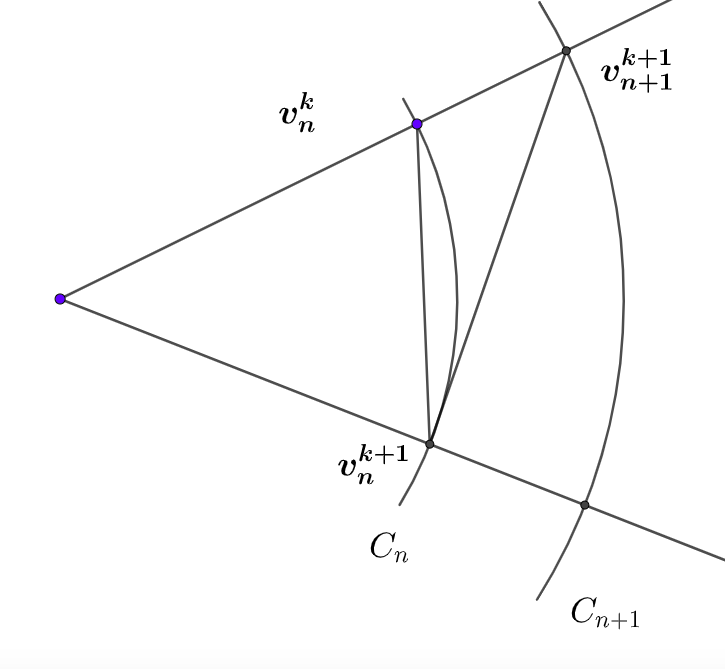}
      \includegraphics[width=0.32\linewidth]{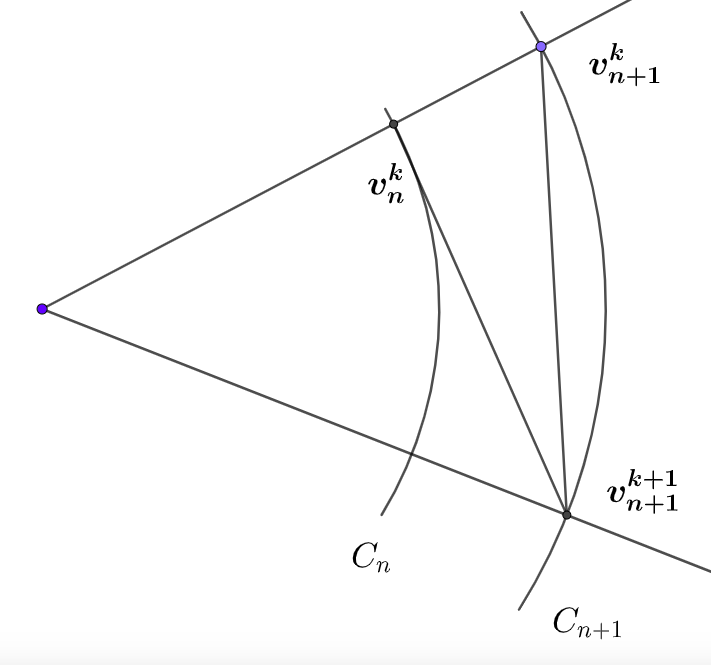}
    \includegraphics[width=0.32\linewidth]{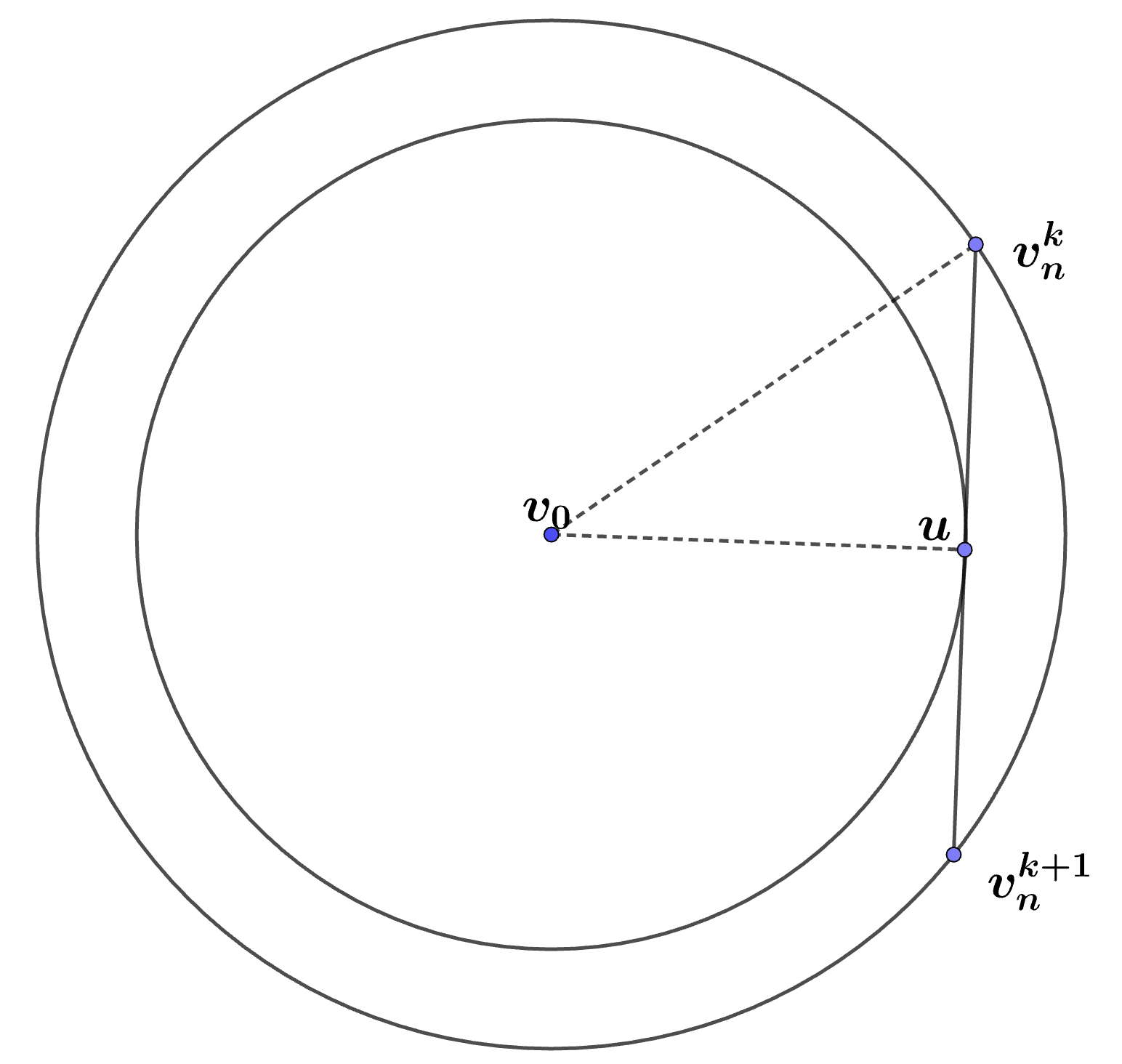}
  \caption{Critical cases for two types of triangles and the longest edge in $A_n$.}
  \label{caseanalysis}
\end{figure}
	In the case for triangle of type 1, we have 
	$$\tanh r_n = \tanh r_{n+1} \cos \frac{\pi}{3n},$$
	and in the case for triangle of type 2, we have 
	$$\tanh r_n = \tanh r_{n+1} \cos \frac{\pi}{3(n+1)}.$$
    To make sure the edges of type 2 stay in $A_n$, we need 
	$$\tanh r_n < \tanh r_{n+1} \cos \frac{\pi}{3n},$$
	which is equivalent to 
	$$\sinh(r_{n+1} - r_n) > \cosh r_n\sinh r_{n+1} (1 - \cos \frac{\pi}{3n}).$$
	One can show that the left side has order $o(1/n)$ and right side has order $o(1/n^{2-2\alpha})$, so the inequality holds for large $n$ if $0<\alpha<1/2$.
	
	This means that there exists $N>0$ such that for $n>N$, the triangles of type 2 will not intersect with the previous layers. 
	For  $0\leq n\leq N$, we can treat it as a triangulated polygon in the Euclidean plane and construct a geometric triangulation arbitrarily. 
\end{proof}

 \begin{lemma}
 The lengths of edges of the geometric triangulation above are uniformly bounded.
 \end{lemma}
\begin{proof}
The length of edges of type 0 is clearly bounded by $r_{n+1}-r_n \leq 1$. The length of the edge of type 1 is bounded by the length of the arc of the circle connecting $v_n^k$ and $v_n^{k+1}$ given by 
$$\sinh l(v_n^kv_n^{k+1}) < \frac{\pi}{3n}\sinh r_n \sim o(\frac{1}{n^{1-\alpha}}) \to 0.$$
From the previous lemma, all the edges of type 2 in the $n$-th layer lie in the ring $A_n$. The length of the edges in $A_n$ with vertices on $C_n$ is bounded by the critical case on the right picture of Figure \ref{caseanalysis}, where the edge $v_n^{k+1}v_n^k$ is tangent to $C_n$. 

The length $L$ of this edge is given by the hyperbolic Pythagorean theorem
$$L  = 2\cosh^{-1}(\frac{\cosh r_{n+1}}{\cosh r_n}) \leq 2\cosh^{-1}(2(1 + \frac{1}{n+1})^\alpha) < \infty.$$
Hence, we can bound the edges of all the three types for large $n>N$. The finite part of the triangulation is naturally bounded. This completes the proof.
\end{proof}

\section{Delaunay triangulations in $\mathbb{H}^2$}

In this section, we prove Theorem \ref{main2} by a simple explicit construction of a family of Delaunay triangulations whose stretch factors tend to infinity, in contrast to the uniform bound on Delaunay triangulations in the Euclidean plane. 

Given a finite point set $V\subset \mathbb{E}^2$, its Delaunay triangulation $T$ can be regarded as a metric graph $(T, d)$ such that the weight on each edge is its length in $\mathbb{E}^2$, and the distance between two vertices of the graph is the shortest path connecting them in $(T, d)$. The \textit{stretch factor} is defined to be
$$s(T) = \max_{p, q\in V}\frac{d(p, q)}{|p-q|},$$
where $|p-q|$ is the Euclidean distance between two distinct $p$ and $q$. It is straightforward to generalize this definition to $\mathbb{H}^2$ by replacing $|p-q|$ with the hyperbolic distance between $p$ and $q$. It is shown \cite{Xia} that $s(T)\leq 1.998$ for any such $V$ in $\mathbb{E}^2$. We construct examples below to show that there is no similar uniform bound in $\mathbb{H}^2$.


\begin{figure}[ht!]
\centering
\begin{overpic}[width=0.6\linewidth]{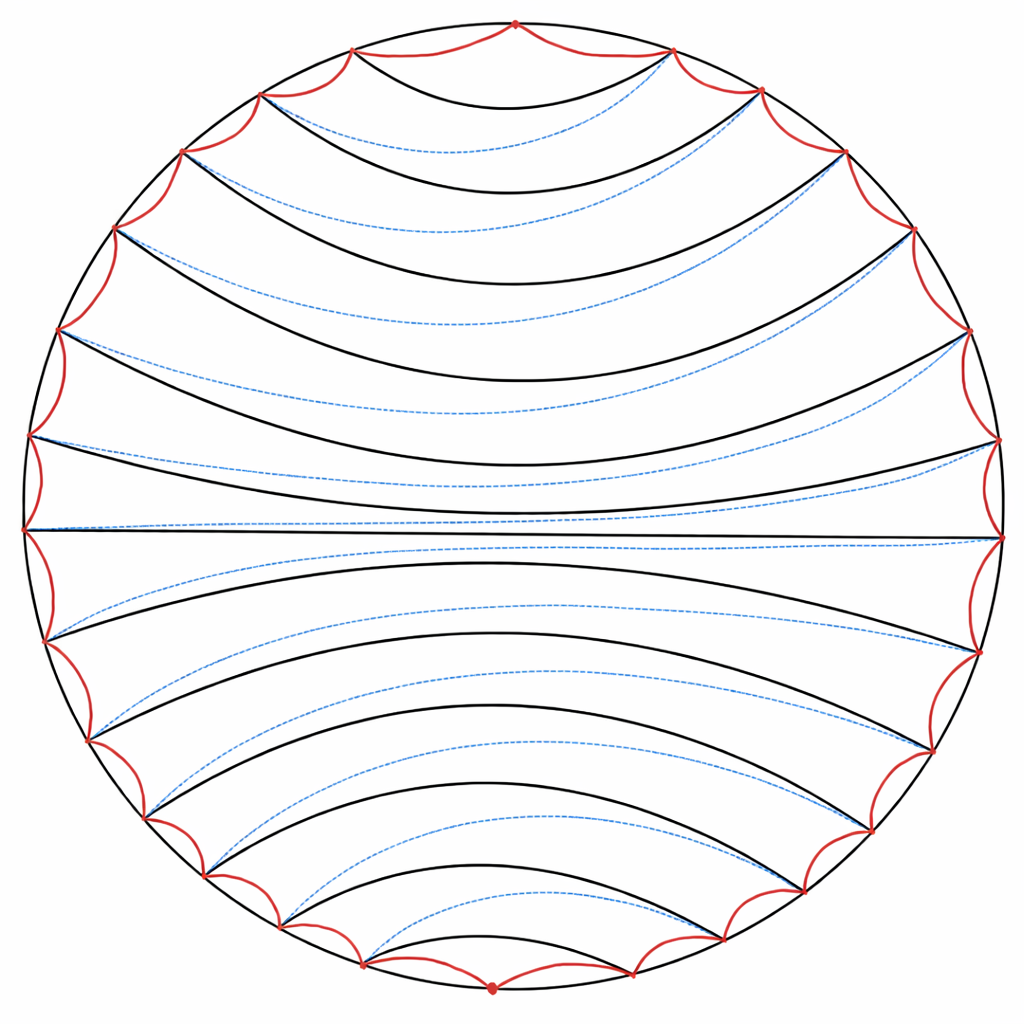}
    \put(48,0){$v_0$}
    \put(60,1){$v_1^R$}
    \put(70,4){$v_2^R$}
    \put(35,1){$v_1^L$}
    \put(24,4){$v_2^L$}
    \put(3,25){$v_k^L$}
    \put(92,25){$v_k^R$}
    \put(-3,48){$v_{n}^L$}
    \put(100,48){$v_{n}^R$}
    \put(48,100){$v_{2n}$}
\end{overpic}
\caption{A Delaunay triangulation of a regular $4n$-gon in $\mathbb{H}^2$.}
  \label{hyperbolicde}
\end{figure}

Let $C$ be a hyperbolic circle centered at the origin of the Poincare disk model of $\mathbb{H}^2$ with radius $\rho$. We construct a regular polygon $P_{4n}$ with $4n$ vertices $V$ evenly distributed on $C$ symmetrically about the $y$-axis.  We denote the vertex at the south pole of the circle as $v_0$ and the vertex at the north pole as $v_{2n}$. Meanwhile, the vertices on the left semicircle are indexed by $v_1^L, v_2^L, \cdots, v_{2n-1}^L$ and the vertices on the right semicircle are indexed by  $v_1^R, v_2^R, \cdots, v_{2n-1}^R$. 
We triangulate $P_{4n}$ by adding its diagonals connecting vertices $v_i^L$ with $v_i^R$ for $i= 1, \cdots, 2n-1$ and vertices $v_i^L$ to $v_{i+1}^R$ for $i = 1, \cdots, 2n-2$. See Figure \ref{hyperbolicde} with $n = 8$. 

Let $T$ denote this triangulation, which is a Delaunay triangulation in $\mathbb{H}^2$ determined by $V$. Indeed, Delaunay triangulations are characterized by the \textit{empty disk property}: the interiors of circumdisks of all triangles contain no points in $V$. The triangulation $T$ is the critical case of Delaunay triangulation where all vertices are on the same circle. 

Let $d$ be the metric on the one-skeleton $T$ induced from edge lengths of $T$ in $\mathbb{H}^2$. Let $l(uv)$ denote the length between two vertices $u$ and $v$ of $P_{4n}$. It is straightforward to check that the $$\lim_{n\to \infty}\sum_{i = 0}^{2n-1}l(v_i^Lv_{i+1}^L) = \lim_{n\to \infty}\sum_{i = 0}^{2n-1}l(v_i^Rv_{i+1}^R) = \pi\sinh \rho.$$
Namely, the perimeters of the regular $4n$-gons approach to the perimeter of the circle. Indeed, a simple computation shows that the boundary edge has length 
$$l(v_i^Lv_{i+1}^L) = 2\sinh^{-1}(\sinh\rho\sin\frac{\pi}{2n}).$$
Clearly, all the diagonals are longer than the boundary length $l(v_i^Lv_{i+1}^L)$. The key observation about $T$ is 
\begin{proposition}
     Shortest paths in $(T, d)$ connecting $v_0$ and $v_{2n}$ contain no diagonal of $T$.
\end{proposition}
\begin{proof}
    Let $\gamma = v_0, \cdots, v_{2n}$ be the shortest path connecting them in $(T, d)$. First, it is easy to see that the list of vertices contain all the indices $0, 1, 2, \cdots, 2n$, possibly with repetition. 
    
    We claim that the indices of the vertices in this path will not decrease. Namely, if $v_k^L$ (or $v_k^R$) is in this path, then the index of the next vertex is at least $k$. Indeed,  assume the opposite and $k\geq 1$, and let $v_k^L$ be a vertex in $\gamma$ such that its consecutive vertex has a smaller index, which must be $v_{k-1}^L$.  There is a copy of $v_k^R$ after this $v_{k-1}^L$ with index $k$ since $\gamma$ is a shortest path, namely,  $$\gamma = v_0\cdots v_k^Lv_{k-1}^L\cdots v_k^R\cdots v_{2n}.$$
    Then we can delete all the vertices between $v_k^L$ and $v_k^R$ to make $\gamma$ shorter, since the edge $v_k^Lv_k^R$ is the shortest geodesic in $\mathbb{H}^2$ connecting them.  The case of $v_k^R$ is similar. 

    Next, we claim that the index in $\gamma$ is strictly increasing, namely, there is no edge $v_k^Lv_k^R$ in $\gamma$. If $v_k^Lv_k^R$ is part of $\gamma$, the previous vertex is $v_{k-1}^L$, since it can not be $v_k^R$ or vertices of larger indices. Then we can replace $v_{k-1}^Lv_k^Lv_k^R$ by $v_{k-1}^Lv_k^R$ to make $\gamma$ shorter.  If $v_k^Rv_k^L$ is part of $\gamma$, we can assume it is the last edge of this type of diagonal in $\gamma$. There are two possibilities: either all the vertices of $\gamma$ after $v_k^R$ are on the left semicircle, or there exists a unique edge $v_m^Lv_{m+1}^R$ in $\gamma$. Indeed, since $v_k^Rv_k^L$ is the last edge of this type, $\gamma$ can not travel from right to left after it, and  
    $$\gamma = v_0\cdots v_k^Rv_k^Lv_{k+1}^L\cdots v_m^Lv_{m+1}^Rv_{m+2}^R\cdots v_{2n}.$$
    However, we can construct a shorter path by replacing $v_k^Lv_{k+1}^L\cdots v_m^L$ by $v_{k+1}^R\cdots v_m^R$ in $\gamma$, since 
    $$\sum_{i = k}^{m-1}l(v_i^Rv_{i+1}^R) = \sum_{i = k}^{m-1}l(v_i^Lv_{i+1}^L), \quad l(v_m^Lv_{m+1}^R)> l(v_m^Rv_{m+1}^R).$$
    
    Finally, we show that there is no diagonal $v_k^Lv_{k+1}^R$ in $\gamma$. Notice that if $\gamma$ contains such an edge, it is unique, since $\gamma$ contains no edge traveling from right to left. Assume this case, then 
    $$\gamma = v_0v_1^L\cdots v_k^Lv_{k+1}^Rv_{k+2}^R\cdots v_{2n}.$$
    However, it is longer than the path $\gamma' = v_0v_1^R\cdots v_k^Rv_{k+1}^Rv_{k+2}^R\cdots v_{2n}$, since by symmetry, the only difference between their lengths is given by $l(v_k^Lv_{k+1}^R) >l(v_k^Rv_{k+1}^R)$. This shows that $\gamma'$ or its symmetric path on the left semicircle is the shortest path between $v_0$ and $v_{2n}$. 
\end{proof}
\begin{theorem}
    For any $C>0$, there is a Delaunay triangulation $T$ in $\mathbb{H}^2$ such that the stretch factor $s(T)>C$.
\end{theorem}
\begin{proof}
    Consider the triangulation $T$ of $P_{4n}$ above and choose $\rho$ large enough such that 
    $$\frac{\pi\sinh \rho}{4\rho}>C.$$
     For this $\rho$, we can choose $n$ sufficiently large such that the shortest path connecting $v_0$ and $v_{2n}$ in $T$ has length larger than $(\pi\sinh\rho)/2$.
     Since the hyperbolic distance between $v_0$ and $v_{2n}$ in $\mathbb{H}^2$ is $2\rho$, we have 
    $$ s(T)\geq\frac{d(v_0, v_{2n})}{2\rho} \geq \frac{\pi\sinh \rho}{4\rho}>C. $$
\end{proof}

\end{document}